\documentstyle[12pt]{article}
\date{}

\begin{document}

\title{Robust Performance of A Class of Control Systems}
\author{Long Wang\thanks{Supported by National Natural Science Foundation of China (69925307), National Key Project of China, National Key Basic Research Special Fund (No. G1998020302) and the Research Fund for the Doctoral Program of Higher Education. }\\\\{\small Center for Systems and Control, Department of Mechanics and Engineering Science,}\\
{\small Peking University, Beijing 100871, CHINA}}

\maketitle
%\baselineskip 24pt
%\large

\begin{abstract}

Some Kharitonov-like robust Hurwitz stability criteria are established for a class of complex polynomial families with nonlinearly correlated perturbations. These results are extended to the polynomial matrix case and non-interval D-stability case. Applications of these results in testing of robust strict positive realness of real and complex interval transfer function families are also presented.

Keywords: Uncertain Systems, Robustness Analysis, Kharitonov's Theorem, Complex Interval Polynomials, Polynomial Matrix Family, Hurwitz Stability, D-Stability, Transfer Functions, Strict Positive Realness.
\end{abstract}

\vspace{3mm}
\section{Introduction}

\hspace{1cm} Motivated by the seminal theorem of Kharitonov on robust stability of interval 
polynomials\cite{Khar, Khar1}, a number of papers on robustness analysis of uncertain
systems have been published in the past few years\cite{Holl, Bart, Fu, wang1, 
wang2, Barmish, Chap, wang3}. Kharitonov's theorem states that the Hurwitz stability of
the real (or complex) interval polynomial family can be guaranteed by the Hurwitz 
stability of four (or eight) prescribed critical vertex polynomials in this family.
This result is significant since it reduces checking stability of infinitely many 
polynomials to checking stability of finitely many polynomials, and the number of
critical vertex polynomials need to be checked is independent of the order of the 
polynomial family. An important extension of Kharitonov's theorem is the edge theorem
discovered by Bartlett, Hollot and Huang\cite{Bart}. The edge theorem states that
the stability of a polytope of polynomials can be guaranteed by the stability of its
one-dimensional exposed edge polynomials. The significance of the edge theorem is 
that it allows some (affine) dependency among polynomial coefficients, and applies to
more general stability regions, e.g., unit circle, left sector, shifted half plane,
hyperbola region, etc. When the dependency among polynomial coefficients is nonlinear,
however, Ackermann shows that checking a subset of a polynomial family generally can 
not guarantee the stability of the entire family\cite{Ack1, Ack2}.

In this paper, we consider a class of complex polynomial families with nonlinear 
coefficient dependency. Based on our previous results, we will establish some 
Kharitonov-like robust stability criteria, i.e., the entire family is stable if 
and only if some critical vertices in this family are stable, and the number of 
critical vertices is independent of the order of the polynomial family. We will then 
extend our results to the polynomial matrix case and non-interval D-stability case.
Applications of these results in testing strict positive realness of interval transfer 
function family are also presented.

\vspace{3mm}
\section{Main Results}

A polynomial $p(s)$ is said to be Hurwitz stable, denoted by $p(s) \in H$, if all its roots lie within the open left half of the complex plane ${\bf C}$. A polynomial family $P$ is said to be Hurwitz stable, denoted by $P \subset H$, if all polynomials in $P$ are Hurwitz stable.

Consider the $n$-th order real interval polynomial family

\begin{equation}
\Gamma = \left \{ p(s) \mid p(s) = \sum_{i=0}^{n} q_{i} s^{i}\;,\;q_{i} \in [q_{i}^{-}\;,\;q_{i}^{+}]\;,\;i=0,1,\cdots ,n \right \}
\end{equation}

\noindent and define the four Kharitonov polynomials of $\Gamma$ as

\begin{equation}
K_1 (s) = q_0^- + q_1^- s + q_2^+ s^2 + q_3^+ s^3 + q_4^- s^4 + q_5^- s^5 + \cdots
\end{equation}

\begin{equation}
K_2 (s) = q_0^+ + q_1^+ s + q_2^- s^2 + q_3^- s^3 + q_4^+ s^4 + q_5^+ s^5 + \cdots
\end{equation}

\begin{equation}
K_3 (s) = q_0^+ + q_1^- s + q_2^- s^2 + q_3^+ s^3 + q_4^+ s^4 + q_5^- s^5 + \cdots
\end{equation}

\begin{equation}
K_4 (s) = q_0^- + q_1^+ s + q_2^+ s^2 + q_3^- s^3 + q_4^- s^4 + q_5^+ s^5 + \cdots
\end{equation}

\vspace{0.6cm}

\noindent {\bf Lemma 1}(Kharitonov's Theorem for Real Polynomials)\cite{Khar}

\begin{equation}
\Gamma \subset H \Longleftrightarrow K_1 (s)\;,\;K_2 (s)\;,\;K_3 (s)\;,\;K_4 (s) \in H
\end{equation}

Consider the $n$-th order complex interval polynomial family

\begin{equation}
\Delta = \left \{ \delta (s) \mid \delta (s) = \sum_{i=0}^{n} (\alpha_{i} + j \beta_{i}) s^{i}\;,\;\alpha_{i} \in [\alpha_{i}^{-}\;,\;\alpha_{i}^{+}]\;,\;\;\beta_{i} \in [\beta_{i}^{-}\;,\;\beta_{i}^{+}]\;,\;i=0,1,\cdots ,n \right \}
\end{equation}

\noindent and define the eight Kharitonov polynomials of $\Delta$ as

{\small

\begin{equation}
K_1^+ (s) = (\alpha^-_0 + j \beta^-_0) + (\alpha^-_1 + j \beta^+_1) s + (\alpha^+_2 + j \beta^+_2) s^2 + (\alpha^+_3 + j \beta^-_3) s^3 + (\alpha^-_4 + j \beta^-_4) s^4 + (\alpha^-_5 + j \beta^+_5) s^5 + \cdots
\end{equation}

\begin{equation}
K_2^+ (s) = (\alpha^-_0 + j \beta^+_0) + (\alpha^+_1 + j \beta^+_1) s + (\alpha^+_2 + j \beta^-_2) s^2 + (\alpha^-_3 + j \beta^-_3) s^3 + (\alpha^-_4 + j \beta^+_4) s^4 + (\alpha^+_5 + j \beta^+_5) s^5 + \cdots
\end{equation}

\begin{equation}
K_3^+ (s) = (\alpha^+_0 + j \beta^-_0) + (\alpha^-_1 + j \beta^-_1) s + (\alpha^-_2 + j \beta^+_2) s^2 + (\alpha^+_3 + j \beta^+_3) s^3 + (\alpha^+_4 + j \beta^-_4) s^4 + (\alpha^-_5 + j \beta^-_5) s^5 + \cdots
\end{equation}

\begin{equation}
K_4^+ (s) = (\alpha^+_0 + j \beta^+_0) + (\alpha^+_1 + j \beta^-_1) s + (\alpha^-_2 + j \beta^-_2) s^2 + (\alpha^-_3 + j \beta^+_3) s^3 + (\alpha^+_4 + j \beta^+_4) s^4 + (\alpha^+_5 + j \beta^-_5) s^5 + \cdots
\end{equation}

\begin{equation}
K_1^- (s) = (\alpha^-_0 + j \beta^-_0) + (\alpha^+_1 + j \beta^-_1) s + (\alpha^+_2 + j \beta^+_2) s^2 + (\alpha^-_3 + j \beta^+_3) s^3 + (\alpha^-_4 + j \beta^-_4) s^4 + (\alpha^+_5 + j \beta^-_5) s^5 + \cdots
\end{equation}

\begin{equation}
K_2^- (s) = (\alpha^-_0 + j \beta^+_0) + (\alpha^-_1 + j \beta^-_1) s + (\alpha^+_2 + j \beta^-_2) s^2 + (\alpha^+_3 + j \beta^+_3) s^3 + (\alpha^-_4 + j \beta^+_4) s^4 + (\alpha^-_5 + j \beta^-_5) s^5 + \cdots
\end{equation}

\begin{equation}
K_3^- (s) = (\alpha^+_0 + j \beta^-_0) + (\alpha^+_1 + j \beta^+_1) s + (\alpha^-_2 + j \beta^+_2) s^2 + (\alpha^-_3 + j \beta^-_3) s^3 + (\alpha^+_4 + j \beta^-_4) s^4 + (\alpha^+_5 + j \beta^+_5) s^5 + \cdots
\end{equation}

\begin{equation}
K_4^- (s) = (\alpha^+_0 + j \beta^+_0) + (\alpha^-_1 + j \beta^+_1) s + (\alpha^-_2 + j \beta^-_2) s^2 + (\alpha^+_3 + j \beta^-_3) s^3 + (\alpha^+_4 + j \beta^+_4) s^4 + (\alpha^-_5 + j \beta^+_5) s^5 + \cdots
\end{equation}

}

\vspace{0.6cm}

\noindent {\bf Lemma 2}(Kharitonov's Theorem for Complex Polynomials)\cite{Khar1}

\begin{equation}
\Delta \subset H \Longleftrightarrow K_1^+ (s)\;,\;K_2^+ (s)\;,\;K_3^+ (s)\;,\;K_4^+ (s)\;,\;K_1^- (s)\;,\;K_2^- (s)\;,\;K_3^- (s)\;,\;K_4^- (s) \in H
\end{equation}

Now consider the $n_u$-th\,,\,$n_v$-th order real interval polynomial families $\Gamma_u$ and $\Gamma_v$\,. Denote their Kharitonov polynomials as $K_i^u (s)\;,\;i=1,2,3,4$ and $K_j^v (s)\;,\;j=1,2,3,4$ respectively.

Similarly, consider the $n_u$-th\,,\,$n_v$-th order complex interval polynomial families $\Delta_u$ and $\Delta_v$\,. Denote their Kharitonov polynomials as $K_i^{+u} (s)\;,\;K_i^{-u} (s)\;,\;i=1,2,3,4$ and $K_j^{+v} (s)\;,\;K_j^{-v} (s)\;,\;j=1,2,3,4$ respectively.

For any function $f(x,y)$, define

\begin{equation}
  f(\Gamma_u\,,\,\Gamma_v) = \{ f(p_u (s)\,,\,p_v (s)) \mid p_u (s) \in \Gamma_u\,,\,p_v (s) \in \Gamma_v \}
\end{equation}

\begin{equation}
  f(\Delta_u\,,\,\Delta_v) = \{ f(\delta_u (s)\,,\,\delta_v (s)) \mid \delta_u (s) \in \Delta_u\,,\,\delta_v (s) \in \Delta_v \}
\end{equation}

\vspace{0.6cm}

\noindent {\bf Lemma 3}\cite{Holl}

For any fixed complex number $z \in {\bf C}$, suppose the polynomial family $\Gamma_u - z \Gamma_v$ has a fixed order. Then

\begin{equation}
  \begin{array}{c}
  \Gamma_u -z \Gamma_v \subset H \Longleftrightarrow\\[6mm]
  K_i^u (s) - z K_j^v (s) \in H\,,\,\,\,\,\,\,\,\,i,j=1,2,3,4
  \end{array}
\end{equation}

If the location of $z$ is known, then the number of critical vertices need to be checked can further be reduced. For example, if $z$ is on the negative real axis, then only 4 out of the 16 critical vertices need to be checked, namely

\begin{equation}
  \begin{array}{c}
  \Gamma_u -z \Gamma_v \subset H \Longleftrightarrow\\[6mm]
  K_i^u (s) - z K_i^v (s) \in H\,,\,\,\,\,\,\,\,\,i=1,2,3,4;
  \end{array}
\end{equation}

\noindent if $z$ is on the imaginary axis, then only 8 critical vertices need to be checked; if $z$ is in the left half of the complex plane, then only 12 critical vertices need to be checked\cite{Holl, wang1, wang2, Barmish}.

For complex polynomials, we have the following similar result

\noindent {\bf Lemma 4}

For any fixed complex number $z \in {\bf C}$, suppose the polynomial family $\Delta_u - z \Delta_v$ has a fixed order. Then

\begin{equation}
  \begin{array}{c}
  \Delta_u -z \Delta_v \subset H \Longleftrightarrow\\[6mm]
  K_i^{+u} (s) - z K_j^{+v} (s),\,\,\,\,  K_i^{-u} (s) - z K_j^{-v} (s) \in H\,,\,\,\,\,\,\,\,\,i,j=1,2,3,4
  \end{array}
\end{equation}

\vspace{0.6cm}

\noindent {\bf Theorem 1}

Consider the polynomial family

\begin{equation}
  a_m \Gamma_u^m + a_{m-1} \Gamma_u^{m-1} \Gamma_v + a_{m-2} \Gamma_u^{m-2} \Gamma_v^2 + \cdots \cdots + a_2 \Gamma_u^2 \Gamma_v^{m-2} + a_1 \Gamma_u \Gamma_v^{m-1} + a_0 \Gamma_v^m
\end{equation}

\noindent where $a_k \in {\bf R}\,,\,k = 0,1, \cdots ,m\,.$ Suppose it has a fixed order. Then

\begin{equation}
  \begin{array}{l}
  a_m \Gamma_u^m + a_{m-1} \Gamma_u^{m-1} \Gamma_v + \cdots \cdots + a_1 \Gamma_u \Gamma_v^{m-1} + a_0 \Gamma_v^m \subset H \Longleftrightarrow\\[6mm]
  a_m [K_i^u (s)]^m + a_{m-1} [K_i^u (s)]^{m-1} K_j^v (s) + \cdots \cdots\\[6mm]
 \,\,\,\,\,\,\,\,\,\,\,\,\,\,\,\,\,\,\,\, + a_1 K_i^u (s) [K_j^v (s)]^{m-1} + a_0 [K_j^v (s)]^m \in H\,,\,\,\,\,\,\,\,\,i,j=1,2,3,4 
  \end{array}
\end{equation}

\vspace{0.6cm}

\noindent Proof: Consider the polynomial

\begin{equation}
  q(z) = a_m z^m + a_{m-1} z^{m-1} + a_{m-2} z^{m-2} + \cdots \cdots + a_2 z^2 + a_1 z + a_0
\end{equation}

\noindent Let $r = \max \{k \mid a_k \neq 0 \}$. Then $q(z)$ can be expressed as

\begin{equation}
  q(z) = a_r (z-z_1)(z-z_2) \cdots \cdots (z-z_{r-1})(z-z_r)
\end{equation}

\noindent where $z_1\,,\,z_2\,,\, \cdots \cdots \,,\,z_{r-1}\,,\,z_r \in {\bf C}$. Hence, we have 

\begin{equation}
  \begin{array}{cl}
& a_m \Gamma_u^m + a_{m-1} \Gamma_u^{m-1} \Gamma_v + \cdots \cdots + a_1 \Gamma_u \Gamma_v^{m-1} + a_0 \Gamma_v^m \subset H\\[6mm]
\Longleftrightarrow & \Gamma_v^m \left [ a_r \left ( \frac{\Gamma_u}{\Gamma_v} \right )^r + a_{r-1} \left ( \frac{\Gamma_u}{\Gamma_v} \right )^{r-1} + \cdots \cdots + a_1 \left ( \frac{\Gamma_u}{\Gamma_v} \right ) + a_0 \right ] \subset H\\[6mm]
\Longleftrightarrow & \Gamma_v^m \left [ a_r \left ( \frac{\Gamma_u}{\Gamma_v} - z_1 \right ) \left ( \frac{\Gamma_u}{\Gamma_v} - z_2 \right ) \cdots \cdots \left ( \frac{\Gamma_u}{\Gamma_v} - z_{r-1} \right ) \left ( \frac{\Gamma_u}{\Gamma_v} - z_r \right ) \right ] \subset H\\[6mm]
\Longleftrightarrow & a_r \Gamma_v^{m-r} ( \Gamma_u - z_1 \Gamma_v )( \Gamma_u - z_2 \Gamma_v ) \cdots \cdots ( \Gamma_u - z_{r-1} \Gamma_v )( \Gamma_u - z_r \Gamma_v ) \subset H\\[6mm]
\Longleftrightarrow & \left \{ \begin{array}{lr}        \Gamma_u - z_k \Gamma_v \subset H\,,\,k=1,2, \cdots \cdots ,r-1,r & \,\,\,\,\,\,\,\,\,\,\,\,\,\,\,r=m\\[6mm]      \Gamma_u - z_k \Gamma_v \subset H\,,\,k=1,2, \cdots \cdots ,r-1,r\,\,\,\, {\rm and}\,\,\,\, \Gamma_v \subset H   & \,\,\,\,\,\,\,\,\,\,\,\,\,\,\,r<m        \end{array} \right.\\[12mm]
\stackrel{\tiny {\rm Lemmas}\,\,1\&3}{\Longleftrightarrow} & {\scriptsize \left \{ \begin{array}{lr}        K_i^u (s) - z_k K_j^v (s) \in H\,,\,i,j=1,2,3,4\,,\,k=1,2, \cdots ,r-1,r & \,\,\,\,\,\,\,\,\,\,\,\,\,r=m\\[6mm]      K_i^u (s) - z_k K_j^v (s) \in H\,,\,i,j=1,2,3,4\,,\,k=1,2, \cdots ,r-1,r\,\,\,\, {\rm and}\,\,\,\, K_j^v (s) \in H\,,\,j=1,2,3,4   & \,\,\,\,\,\,\,\,\,\,\,\,\,r<m        \end{array} \right.}\\[12mm]
\end{array}
\end{equation}

\begin{equation}
  \begin{array}{cl}
\Longleftrightarrow & a_r [K_j^v (s)]^{m-r} [K_i^u (s) - z_1 K_j^v (s)][K_i^u (s) - z_2 K_j^v (s)] \cdots\\[6mm]
& \,\,\,\,\,\,\,\,\,\,\,\,\,\,\,\,\,\,\,\,\cdots [K_i^u (s) - z_{r-1} K_j^v (s)][K_i^u (s) - z_r K_j^v (s)] \in H\,,\,\,\,\,\,\,\,\,i,j=1,2,3,4\\[6mm]
\Longleftrightarrow & [K_j^v (s)]^m \left [ a_r \left ( \frac{K_i^u (s)}{K_j^v (s)} - z_1 \right ) \left ( \frac{K_i^u (s)}{K_j^v (s)} - z_2 \right ) \cdots \right.\\[6mm]
& \left. \,\,\,\,\,\,\,\,\,\,\,\,\,\,\,\,\,\,\,\,\cdots \left ( \frac{K_i^u (s)}{K_j^v (s)} - z_{r-1} \right ) \left ( \frac{K_i^u (s)}{K_j^v (s)} - z_r \right ) \right ] \in H\,,\,\,\,\,\,\,\,\,i,j=1,2,3,4\\[6mm]
\Longleftrightarrow & [K_j^v (s)]^m \left [ a_r \left ( \frac{K_i^u (s)}{K_j^v (s)} \right )^r + a_{r-1} \left ( \frac{K_i^u (s)}{K_j^v (s)} \right )^{r-1} + \cdots \right.\\[6mm]
& \left. \,\,\,\,\,\,\,\,\,\,\,\,\,\,\,\,\,\,\,\,\cdots + a_1 \left ( \frac{K_i^u (s)}{K_j^v (s)} \right ) + a_0 \right ] \in H\,,\,\,\,\,\,\,\,\,i,j=1,2,3,4\\[6mm]
\Longleftrightarrow & a_m [K_i^u (s)]^m + a_{m-1} [K_i^u (s)]^{m-1} K_j^v (s) + \cdots\\[6mm]
& \,\,\,\,\,\,\,\,\,\,\,\,\,\,\,\,\,\,\,\,\cdots + a_1 K_i^u (s) [K_j^v (s)]^{m-1} + a_0 [K_j^v (s)]^m \in H\,,\,\,\,\,\,\,\,\,i,j=1,2,3,4 
  \end{array}
\end{equation}

\noindent This completes the proof.

\vspace{0.6cm}

From the proof of Theorem 1 and by Lemmas 2 and 4, we have

\vspace{0.3cm}

\noindent {\bf Theorem 2}

Consider the polynomial family

\begin{equation}
  c_m \Delta_u^m + c_{m-1} \Delta_u^{m-1} \Delta_v + c_{m-2} \Delta_u^{m-2} \Delta_v^2 + \cdots \cdots + c_2 \Delta_u^2 \Delta_v^{m-2} + c_1 \Delta_u \Delta_v^{m-1} + c_0 \Delta_v^m
\end{equation}

\noindent where $c_k \in {\bf C}\,,\,k = 0,1, \cdots ,m\,.$ Suppose it has a fixed order. Then

\begin{equation}
  \begin{array}{l}
  c_m \Delta_u^m + c_{m-1} \Delta_u^{m-1} \Delta_v + \cdots \cdots + c_1 \Delta_u \Delta_v^{m-1} + c_0 \Delta_v^m \subset H \Longleftrightarrow\\[6mm]
  c_m [K_i^{+u} (s)]^m + c_{m-1} [K_i^{+u} (s)]^{m-1} K_j^{+v} (s) + \cdots \cdots\\[6mm]
 \,\,\,\,\,\,\,\,\,\,\,\,\,\,\,\,\,\,\,\, + c_1 K_i^{+u} (s) [K_j^{+v} (s)]^{m-1} + c_0 [K_j^{+v} (s)]^m \in H\,,\,\,\,\,\,\,\,\,\\[6mm] 
  c_m [K_i^{-u} (s)]^m + c_{m-1} [K_i^{-u} (s)]^{m-1} K_j^{-v} (s) + \cdots \cdots\\[6mm]
 \,\,\,\,\,\,\,\,\,\,\,\,\,\,\,\,\,\,\,\, + c_1 K_i^{-u} (s) [K_j^{-v} (s)]^{m-1} + c_0 [K_j^{-v} (s)]^m \in H\,,\,\,\,\,\,\,\,\,i,j=1,2,3,4
\end{array}
\end{equation}

\noindent Remark. We have established strong Kharitonov-like criteria for the stability of a class of polynomial families with nonlinearly correlated perturbations. The number of critical polynomials need to be checked is independent of the order of the polynomial family.
 
\noindent Example 1

Consider a negative unity feedback system with the forward path as three same blocks in tandem. Each block consists of an interval plant $\frac{N(s)}{D(s)}$ with negative unity feedback. Then, the characteristic polynomial of the closed-loop system is

\begin{equation}
  [N(s)]^3 + [N(s) + D(s)]^3
\end{equation}

\noindent By Theorem 1, we only need to check 16 vertex systems for the stability of the entire uncertain system family. Furthermore, since all the roots of

\begin{equation}
  q(z) = 2z^3 + 3z^2 + 3z + 1
\end{equation}

\noindent lie within the left half of the complex plane, only 12 out of the 16 vertex systems need to be checked to verify robust stability of the entire system family.

\noindent Example 2

Consider a negative unity feedback system with the forward path as a controller and an interval plant $\frac{N(s)}{D(s)}$ in tandem. The controller is simply a gain $k$\,,\, but can be switched among $\{ k_1\,,\,k_2\,,\, \cdots \cdots \,,\,k_m \}$ under different working conditions. Thus, robust stability of the entire system family is tantamount to

\begin{equation}
  [k_1 N(s) + D(s)][k_2 N(s) + D(s)] \cdots \cdots [k_m N(s) + D(s)] \subset H
\end{equation}

\noindent By Theorem 1, we only need to check 16 vertex systems for the stability of the entire uncertain system family. Furthermore, since all the roots of

\begin{equation}
  q(z) = (k_1 z + 1)(k_2 z + 1) \cdots \cdots (k_m z + 1)
\end{equation}

\noindent lie on the real axis, only 8 out of the 16 vertex systems need to be checked. Moreover, if $k_1,k_2, \cdots \cdots ,k_m$ have the same sign, then only 4 out of the 16 vertex systems need to be checked.

\vspace{3mm}
\section{Some Extensions}

\subsection{Extension to Non-Interval D-Stability Case}

Given any stability region $D$ in the complex plane ${\bf C}$, a polynomial $p(s)$ is said to be D-stable, denoted by $p(s) \in D$, if all its roots lie within $D$. A polynomial family $P$ is said to be D-stable, denoted by $P \subset D$, if all polynomials in $P$ are D-stable.

Let the uncertainty bounding set (hyperbox) be

\begin{equation}
  \begin{array}{r}
  Q = \{ q = ( q_1\,,\,q_2\,,\, \cdots \cdots \,,\,q_l )^T \mid q_i \in [ q_i^-\,,\,q_i^+ ]\,,\,\\[6mm]
  i = 1,2, \cdots \cdots ,l \}
  \end{array}
\end{equation}

\noindent and define its one-dimensional edge set as

{\scriptsize
\begin{equation}
  \begin{array}{r}
  Q_E  =  \{ q = ( q_1\,,\,q_2\,,\, \cdots \cdots \,,\,q_l )^T \mid q_k \in [ q_k^-\,,\,q_k^+ ]\,\,\,\,{\rm for\,\,\,\, some}\,\,\,\,\,\,\,\,\,\,\,\,\,\,\,\,\\[6mm]
     k \in \{ 1,2, \cdots \cdots ,l \} \,\,\,\, {\rm and}\,\,\,\, q_i \in \{ q_i^-\,,\,q_i^+ \}\,\,\,\, {\rm for\,\,\,\, all}\,\,\,\, i \neq k \}
  \end{array}      
\end{equation}
}

Consider the $n_1$-th, $n_2$-th order complex polynomials

\begin{equation}
  n(s,q) = \sum_{i=0}^{n_1} c_i (q) s^i
\end{equation}

\begin{equation}
  d(s,q) = \sum_{j=0}^{n_2} b_j (q) s^j
\end{equation}

\noindent where the complex coefficients $c_i (q)\,,\,b_j (q)$ are affine functions of the uncertain parameters $q = ( q_1\,,\,q_2\,,\, \cdots \cdots \,,\,q_l )^T\,,\,$ respectively.

In the sequel, we will suppose that $D^c$ is a connected set. Note that Hurwitz stability and Schur stability are special cases of D-stability.

\vspace{0.3cm}

\noindent {\bf Lemma 5}

For any fixed complex numbers $z_{01}\,,\,z_{02} \in {\bf C}\,,\,$ suppose the polynomial family $\{ z_{01} n(s,q) + z_{02} d(s,q) \mid q \in Q \}$ has a fixed order. Then

\begin{equation}
  \begin{array}{c}
    \{ z_{01} n(s,q) + z_{02} d(s,q) \mid q \in Q \} \subset D \Longleftrightarrow\\[6mm]
    \{ z_{01} n(s,q) + z_{02} d(s,q) \mid q \in Q_E \} \subset D
  \end{array}
\end{equation}

\vspace{0.3cm}

\noindent Proof: Since the coefficients of $z_{01} n(s,q) + z_{02} d(s,q)$ are also affine functions of $q = ( q_1\,,\,q_2\,,\, \cdots \cdots \,,\,q_l )^T\,,\,$ the result follows directly from the Edge Theorem\cite{Bart,Fu}.

\vspace{3mm}

For notational simplicity, define

{\small
\begin{equation}
  \begin{array}{l}
g(s,q) = a_m [n(s,q)]^m + a_{m-1} [n(s,q)]^{m-1} d(s,q)\\[6mm]
\,\,\,\,\,\,\,\,\,\,\,\,\,\,\,\,\,\,\,\,+ a_{m-2} [n(s,q)]^{m-2} [d(s,q)]^2 + \cdots \cdots\\[6mm]
 \,\,\,\,\,\,\,\,\,\,\,\,\,\,\,\,\,\,\,\, + a_2 [n(s,q)]^2 [d(s,q)]^{m-2} + a_1 n(s,q) [d(s,q)]^{m-1}\\[6mm]
\,\,\,\,\,\,\,\,\,\,\,\,\,\,\,\,\,\,\,\, + a_0 [d(s,q)]^m
  \end{array}
\end{equation}
}

\noindent where $a_k \in {\bf C}\,,\,k = 0,1, \cdots \cdots ,m\,.$

\vspace{0.3cm}

\noindent {\bf Theorem 3}

Consider the polynomial family

\begin{equation}
  \{ g(s,q) \mid q \in Q \}
\end{equation}

\noindent Suppose it has a fixed order. Then

\begin{equation}
  \begin{array}{c}
  \{ g(s,q) \mid q \in Q \} \subset D \Longleftrightarrow\\[6mm]
  \{ g(s,q) \mid q \in Q_E \} \subset D
  \end{array}
\end{equation}

\vspace{0.3cm}

\noindent Proof: Consider the polynomial

{\small
\begin{equation}
  q(z) = a_m z^m + a_{m-1} z^{m-1} + a_{m-2} z^{m-2} + \cdots \cdots + a_2 z^2 + a_1 z + a_0
\end{equation}
}

\noindent Let $r = \max \{k \mid a_k \neq 0 \}$. Then $q(z)$ can be expressed as

\begin{equation}
  q(z) = a_r (z-z_1)(z-z_2) \cdots \cdots (z-z_{r-1})(z-z_r)
\end{equation}

\noindent where $z_1\,,\,z_2\,,\, \cdots \cdots \,,\,z_{r-1}\,,\,z_r \in {\bf C}$. Hence, we have 

{\tiny
\begin{equation}
  \begin{array}{cl}
& \{ g(s,q) \mid q \in Q \} \subset D\,\,\,\,\,\, \Longleftrightarrow\,\,\,\,\,\, g(s,q) \in D\,,\,\forall q \in Q\\[6mm]
\Longleftrightarrow & [d(s,q)]^m \left \{ a_m \left [ \frac{n(s,q)}{d(s,q)} \right ]^m + a_{m-1} \left [ \frac{n(s,q)}{d(s,q)} \right ]^{m-1} + \cdots \right.\\[6mm]
& \left. \,\,\,\,\,\,\,\,\,\,\,\,\,\,\,\,\,\,\,\,\cdots + a_1 \left [ \frac{n(s,q)}{d(s,q)} \right ]^1 + a_0 \right \} \in D\,,\,\forall q \in Q\\[6mm]
\Longleftrightarrow & [d(s,q)]^m \left \{ a_r \left [ \frac{n(s,q)}{d(s,q)} - z_1 \right ] \left [ \frac{n(s,q)}{d(s,q)} - z_2 \right ] \cdots \right.\\[6mm]
& \left. \,\,\,\,\,\,\,\,\,\,\,\,\,\,\,\,\,\,\,\,\cdots \left [ \frac{n(s,q)}{d(s,q)} - z_{r-1} \right ] \left [ \frac{n(s,q)}{d(s,q)} - z_r \right ] \right \} \in D\,,\,\forall q \in Q\\[6mm]
\Longleftrightarrow & a_r [d(s,q)]^{m-r} [n(s,q) - z_1 d(s,q)][n(s,q) - z_2 d(s,q)] \cdots\\[6mm]
& \,\,\,\,\,\,\,\,\,\,\,\,\cdots [n(s,q) - z_{r-1} d(s,q)][n(s,q) - z_r d(s,q)] \in D\,,\,\forall q \in Q\\[8mm]
\Longleftrightarrow & {\scriptsize \left \{ \begin{array}{lr}    \begin{array}{l}    n(s,q) - z_k d(s,q) \in D\,,\,\\
        k=1,2, \cdots \cdots ,r-1,r\,,\,\forall q \in Q \end{array} & \,\,\,\,\,\,\,\,\,\,\,\,\,\,\,r=m\\[6mm]   \begin{array}{l}   n(s,q) - z_k d(s,q) \in D\,,\,\\
        k=1,2, \cdots \cdots ,r-1,r\,,\,\forall q \in Q\,\,\,\,\\
        {\rm and}\,\,\,\, d(s,q) \in D\,,\,\forall q \in Q \end{array} & \,\,\,\,\,\,\,\,\,\,\,\,\,\,\,r<m        \end{array} \right.}\\[18mm]
\stackrel{\tiny {\rm Lemma}\,\,5}{\Longleftrightarrow} & {\scriptsize \left \{ \begin{array}{lr}  \begin{array}{l}      n(s,q) - z_k d(s,q) \in D\,,\,\\
        k=1,2, \cdots \cdots ,r-1,r\,,\,\forall q \in Q_E \end{array} & \,\,\,\,\,\,\,\,\,\,\,\,\,\,\,r=m\\[6mm] \begin{array}{l}     n(s,q) - z_k d(s,q) \in D\,,\,\\
        k=1,2, \cdots \cdots ,r-1,r\,,\,\forall q \in Q_E\,\,\,\,\\
        {\rm and}\,\,\,\, d(s,q) \in D\,,\,\forall q \in Q_E  \end{array} & \,\,\,\,\,\,\,\,\,\,\,\,\,\,\,r<m        \end{array} \right.}\\[16mm]
\Longleftrightarrow & a_r [d(s,q)]^{m-r} [n(s,q) - z_1 d(s,q)][n(s,q) - z_2 d(s,q)] \cdots\\[6mm]
& \,\,\,\,\,\,\,\,\,\,\,\,\cdots [n(s,q) - z_{r-1} d(s,q)][n(s,q) - z_r d(s,q)] \in D\,,\,\forall q \in Q_E\\[6mm]
\Longleftrightarrow & [d(s,q)]^m \left \{ a_r \left [ \frac{n(s,q)}{d(s,q)} - z_1 \right ] \left [ \frac{n(s,q)}{d(s,q)} - z_2 \right ] \cdots \right.\\[6mm]
& \left. \,\,\,\,\,\,\,\,\,\,\,\,\,\,\,\,\,\,\,\,\cdots \left [ \frac{n(s,q)}{d(s,q)} - z_{r-1} \right ] \left [ \frac{n(s,q)}{d(s,q)} - z_r \right ] \right \} \in D\,,\,\forall q \in Q_E\\[6mm]
\Longleftrightarrow & [d(s,q)]^m \left \{ a_m \left [ \frac{n(s,q)}{d(s,q)} \right ]^m + a_{m-1} \left [ \frac{n(s,q)}{d(s,q)} \right ]^{m-1} + \cdots \right.\\[6mm]
& \left. \,\,\,\,\,\,\,\,\,\,\,\,\,\,\,\,\,\,\,\,\cdots + a_1 \left [ \frac{n(s,q)}{d(s,q)} \right ]^1 + a_0 \right \} \in D\,,\,\forall q \in Q_E\\[6mm]
\Longleftrightarrow & g(s,q) \in D\,,\,\forall q \in Q_E\,\,\,\,\,\, \Longleftrightarrow\,\,\,\,\,\, \{ g(s,q) \mid q \in Q_E \} \subset D 
  \end{array}
\end{equation}
}

\noindent This completes the proof.

\vspace{6mm}

\noindent Remark. Theorem 3 reveals that, for a class of polynomial family with nonlinearly correlated perturbations, D-stability of the entire family can be ascertained by only checking one-dimensional edge polynomials in this family.

\vspace{6mm}

\subsection{Extension to Polynomial Matrix Families}

Consider the uncertain polynomial matrix

{\tiny
\begin{equation}
  M(s,q) =  \left [ \begin{array}{ccc}
      2n(s,q) + 3d(s,q)\,\,\, &\,\,\, 3n(s,q) + 4d(s,q)\,\,\, &\,\,\, 0\\[6mm]
      0\,\,\,                 &\,\,\, 4n(s,q) + 5d(s,q)\,\,\,  &\,\,\, 2n(s,q)\\[6mm]
      9d(s,q)\,\,\,           &\,\,\, 6n(s,q)\,\,\,           &\,\,\, 5n(s,q) + 6d(s,q)
\end{array} \right ]
\end{equation}
}

\noindent it is easy to see that

%{\tiny
\begin{equation}
  \begin{array}{r}
  \det [M(s,q)] = 16 [n(s,q)]^3 + 176 [n(s,q)]^2 d(s,q)\\[6mm]
  + 279 n(s,q) [d(s,q)]^2 + 90 [d(s,q)]^3
  \end{array}
\end{equation}
%}

\noindent By Theorem 3, we have

\begin{equation}
  \begin{array}{c}
  \{ \det [M(s,q)] \mid q \in Q \} \subset D \Longleftrightarrow\\[6mm]
  \{ \det [M(s,q)] \mid q \in Q_E \} \subset D
  \end{array}
\end{equation}

\noindent Namely, robust D-stability of the entire polynomial matrix family can be ascertained by only checking one-dimensional edges. More generally, for any uncertain polynomial matrix of the form

\begin{equation}
  M(s,q) = [ \alpha_{ij} n(s,q) + \beta_{ij} d(s,q) ]_{n \times n}
\end{equation}

\noindent it is easy to see that the above edge result also holds. Moreover, if $n(s,q)\,,\,d(s,q)$ are replaced by interval polynomial families $\Gamma_u$\,,\,$\Gamma_v$ or $\Delta_u$\,,\,$\Delta_v$ as defined in the last section, then Kharitonov-like results can be established for robust Hurwitz stability of the corresponding polynomial matrix families.

\vspace{0.3cm}

\noindent {\bf Theorem 4}

Consider the polynomial matrix family

\begin{equation}
  M(\delta_u (s), \delta_v (s)) = [ \gamma_{ij} \delta_u (s) + \eta_{ij} \delta_v (s) ]_{n \times n}
\end{equation}

\noindent where $\delta_u (s) \in \Delta_u\,,\,\delta_v (s) \in \Delta_v$, and $\gamma_{ij}, \eta_{ij}, i,j=1,2,......,n$ are complex numbers. Then

\begin{equation}
  \begin{array}{c}
  \{ \det [M(\delta_u (s), \delta_v (s))] \mid \delta_u (s) \in \Delta_u\,,\,\delta_v (s) \in \Delta_v \} \subset H \Longleftrightarrow\\[6mm]
  \{ \det [M(K_i^{+u} (s)\;,\;K_j^{+v} (s))] \mid i,j=1,2,3,4 \} \bigcup\\[6mm]
  \,\,\,\,\,\,\,\,\,\,\,\,\,\,\,\,\,\,\,\,\,\,\,\,\,\,\,\,\,\,\,\,\,\{ \det [M(K_i^{-u} (s)\;,\;K_j^{-v} (s))] \mid i,j=1,2,3,4 \} \subset H
  \end{array}
\end{equation}

\vspace{6mm}
\section{Some Applications}

A proper transfer function $\frac{p(s)}{q(s)}$ is said to be strictly positive real, denoted by $\frac{p(s)}{q(s)} \in SPR$, if

\begin{equation}
  \begin{array}{l}
    1)\,\,q(s) \in H\\[3mm]
    2)\,\,\Re {\frac{p(j \omega)}{q(j \omega)}}>0\;,\;\forall \omega \in R
  \end{array}
\end{equation}

\noindent Suppose $p(s), q(s)$ have positive leading coefficients. Then, it is easy to see that

\begin{equation}
  \frac{p(s)}{q(s)} \in SPR \Longleftrightarrow \lambda p^2 (s) + (1- \lambda) q^2 (s) \in H, \lambda \in [0,1]
\end{equation}

Now consider the proper interval transfer function family

\begin{equation}
  T = \left \{ \frac{p_u (s)}{p_v (s)} \mid p_u (s) \in \Gamma_u\;,\; p_v (s) \in \Gamma_v \right \}
\end{equation}

\noindent In order to have

\begin{equation}
\frac{p_u (s)}{p_v (s)} \in SPR\;,\; p_u (s) \in \Gamma_u\;,\; p_v (s) \in \Gamma_v
\end{equation}

\noindent we must have

\begin{equation}
\lambda \Gamma_u^2 + (1- \lambda) \Gamma_v^2 \subset H, \lambda \in [0,1]
\end{equation}

\noindent Since $\lambda z^2 + (1- \lambda)$ has purely imaginary roots. By Theorem 1, we only need to have\cite{Holl, wang1, wang2, Barmish}

\begin{equation}
\lambda [K^u_1 (s)]^2 + (1- \lambda) [K^v_4 (s)]^2 \in H, \lambda \in [0,1]
\end{equation}

\begin{equation}
\lambda [K^u_2 (s)]^2 + (1- \lambda) [K^v_3 (s)]^2 \in H, \lambda \in [0,1]
\end{equation}

\begin{equation}
\lambda [K^u_3 (s)]^2 + (1- \lambda) [K^v_1 (s)]^2 \in H, \lambda \in [0,1]
\end{equation}

\begin{equation}
\lambda [K^u_4 (s)]^2 + (1- \lambda) [K^v_2 (s)]^2 \in H, \lambda \in [0,1]
\end{equation}

\begin{equation}
\lambda [K^u_1 (s)]^2 + (1- \lambda) [K^v_3 (s)]^2 \in H, \lambda \in [0,1]
\end{equation}

\begin{equation}
\lambda [K^u_2 (s)]^2 + (1- \lambda) [K^v_4 (s)]^2 \in H, \lambda \in [0,1]
\end{equation}

\begin{equation}
\lambda [K^u_3 (s)]^2 + (1- \lambda) [K^v_2 (s)]^2 \in H, \lambda \in [0,1]
\end{equation}

\begin{equation}
\lambda [K^u_4 (s)]^2 + (1- \lambda) [K^v_1 (s)]^2 \in H, \lambda \in [0,1]
\end{equation}

\noindent Equivalently

\begin{equation}
  \begin{array}{l}
 \frac{K_1^u (s)}{K_4^v (s)}\;,\; \frac{K_2^u (s)}{K_3^v (s)}\;,\; \frac{K_3^u (s)}{K_1^v (s)}\;,\; \frac{K_4^u (s)}{K_2^v (s)}\;,\;\\[4mm]
  \;\;\;\;\;\;\; \frac{K_1^u (s)}{K_3^v (s)}\;,\; \frac{K_2^u (s)}{K_4^v (s)}\;,\; \frac{K_3^u (s)}{K_2^v (s)}\;,\; \frac{K_4^u (s)}{K_1^v (s)}\; \in \;SPR
  \end{array}
\end{equation}

\noindent Namely, in order to guarantee that every member of the interval transfer function family $T$ is strictly positive real, we only need to check eight specially selected vertex transfer functions. That is

\begin{equation}
  \begin{array}{l}
  \frac{p_u (s)}{p_v (s)} \in SPR\;,\;\forall p_u (s) \in \Gamma_u\;,\;\forall p_v (s) \in \Gamma_v\\[4mm] 
  \Longleftrightarrow \frac{K_1^u (s)}{K_4^v (s)}\;,\; \frac{K_2^u (s)}{K_3^v (s)}\;,\; \frac{K_3^u (s)}{K_1^v (s)}\;,\; \frac{K_4^u (s)}{K_2^v (s)}\;,\;\\[4mm]
  \;\;\;\;\;\;\;\; \frac{K_1^u (s)}{K_3^v (s)}\;,\; \frac{K_2^u (s)}{K_4^v (s)}\;,\; \frac{K_3^u (s)}{K_2^v (s)}\;,\; \frac{K_4^u (s)}{K_1^v (s)}\; \in \;SPR
  \end{array}
\end{equation}

\noindent which is consistent with the result of Chapellat et al\cite{Chap}, and Wang\cite{wang3}.

Moreover, for any $\gamma \in R\,$, in order to have

\begin{equation}
\gamma + \frac{p_u (s)}{p_v (s)} \in SPR\;,\; p_u (s) \in \Gamma_u\;,\; p_v (s) \in \Gamma_v
\end{equation}

\noindent we must have

\begin{equation}
\lambda [\gamma \Gamma_v + \Gamma_u]^2 + (1- \lambda) \Gamma_v^2 \subset H, \lambda \in [0,1]
\end{equation}

\noindent Since $\lambda (\gamma + z)^2 + (1- \lambda)$ has roots either at first and fourth quadrants (when $\gamma < 0$) or at second and third quadrants (when $\gamma > 0$). By Theorem 1, we only need to check twelve vertices to guarantee robust stability\cite{Holl, wang1, wang2, Barmish}. Namely, in order to guarantee that

\begin{equation}
\gamma + \frac{p_u (s)}{p_v (s)} \in SPR\;,\; p_u (s) \in \Gamma_u\;,\; p_v (s) \in \Gamma_v
\end{equation}

\noindent we only need to check the same property for twelve specially selected vertex transfer functions.

\noindent Remark. The above result can be easily extended to the case of complex interval transfer function family. Namely, every member in the complex interval transfer function family is strictly positive real, if and only if, sixteen specially selected vertex transfer functions in this family are strictly positive real\cite{wang3}.

\vspace{6mm}
\section{Robust Sensitivity Functions}

Denote the $m$-th, $n$-th $(m<n)$ order real interval polynomial families $K_g (s)$, $K_f (s)$ as

\begin{equation}
K_g(s)=\{g(s)|g(s)=\sum_{i=0}^mb_is^i, b_i\in [\underline{b_i},
\overline{b_i}], i=0,1, ......,m\},
\end{equation}

\begin{equation}
K_f(s)=\{f(s)|f(s)=\sum_{i=0}^na_is^i, a_i\in [\underline{a_i},
\overline{a_i}], i=0,1, ......,n\}.
\end{equation}

\noindent For any $f(s)\in K_f(s)$ , it can be expressed as

\begin{equation}
f(s)=\alpha _f(s^2)+s\beta _f(s^2),
\end{equation}

\noindent where

\begin{equation}
\alpha _f(s^2)=a_0+a_2s^2+a_4s^4+a_6s^6+......,
\end{equation}

\begin{equation}
\beta _f(s^2)=a_1+a_3s^2+a_5s^4+a_7s^6+.......
\end{equation}

\noindent Obviously, for any fixed $\omega \in R$, $\alpha _f(-\omega ^2)$ and $\omega \beta _f(-\omega^2)$ are the real and imaginary parts of $f(j\omega )\in C$ respectively.

For the interval polynomial family $K_f (s)$, define

\begin{equation}
\alpha _f^{(1)}(s^2)=\underline{a_0}+\overline{a_2}s^2+\underline{a_4}s^4+%
\overline{a_6}s^6+......,
\end{equation}

\begin{equation}
\alpha _f^{(2)}(s^2)=\overline{a_0}+\underline{a_2}s^2+\overline{a_4}s^4+%
\underline{a_6}s^6+......,
\end{equation}

\begin{equation}
\beta _f^{(1)}(s^2)=\underline{a_1}+\overline{a_3}s^2+\underline{a_5}s^4+%
\overline{a_7}s^6+......,
\end{equation}

\begin{equation}
\beta _f^{(2)}(s^2)=\overline{a_1}+\underline{a_3}s^2+\overline{a_5}s^4+%
\underline{a_7}s^6+......,
\end{equation}

\noindent and denote the four Kharitonov vertex polynomials of $K_f (s)$ as

\begin{equation}
f_{ij}(s)=\alpha _f^{(i)}(s^2)+s\beta _f^{(j)}(s^2), \,\,\,\,\,\,\,\,\,\,\,\, i,j=1,2
\end{equation}

For the interval polynomial family $K_g (s)$, the corresponding $\alpha _g^{(i)} (s), \beta _g^{(j)} (s)$ and $g_{ij}(s) \in K_g (s)$ can be defined analogously.

\vspace{0.4cm}

\noindent {\bf Lemma 6}\cite{das}

For any fixed $\omega \in R$, $f(s)\in K_f(s)$, we have

\begin{equation}
\alpha _f^{(1)}(-\omega ^2)\leq \alpha _f(-\omega ^2)\leq \alpha
_f^{(2)}(-\omega ^2),
\end{equation}

\begin{equation}
\beta _f^{(1)}(-\omega ^2)\leq \beta _f(-\omega ^2)\leq \beta
_f^{(2)}(-\omega ^2).
\end{equation}

\vspace{0.4cm}

\noindent {\bf Lemma 7}\cite{Ack1} (Zero Exclusion Principle)

For the $n$-th order polynomial family

\begin{equation}
f(s,T)=: \{f(s,t)|t\in T\},
\end{equation}

\noindent where $T$ is a bounded connected closed set, and the coefficients of $f(s,t)$ are continuous functions of $t$, then $f(s,T)\in H$ if and only if

1)\,\,\,\,\,there exists $t^{*}\in T,$  such that $f(s,t^{*})\in H$;

2)\,\,\,\,\,$0\notin f(j\omega ,T),$ $\forall \omega \in R$.

Consider the strictly proper open-loop transfer function

\begin{equation}
P=\frac{g(s)}{f(s)}
\end{equation}

\noindent and suppose the closed-loop system is stable under negative unity feedback. Denote its sensitivity function as

\begin{equation}
S=\frac{1}{1+P}=\frac{f(s)}{f(s)+g(s)}
\end{equation}

\noindent Apparently, we have

\begin{equation}
||S||_{\infty} \geq 1
\end{equation}

For notational simplicity, define

\begin{equation}
J_{i_1 j_1 i_2 j_2} (s) = g_{i_1 j_1}(s) + (1 + \delta e^{j \theta}) f_{i_2 j_2}(s), \,\,\,\,\,\,\,\,\,\, \delta \in (0, 1), \,\,\, i_1, j_1, i_2, j_2 = 1,2, \,\,\, \theta \in [-\pi, \pi].
\end{equation}

\vspace{0.4cm}

\noindent {\bf Lemma 8}

Suppose $g(s) + f(s) \in H$. Then, for any $\gamma > 1$, we have

\begin{equation}
||S||_{\infty}<\gamma \Longleftrightarrow g(s) + (1 + \frac{1}{\gamma} e^{j \theta}) f(s) \in H, \,\,\,\,\, \forall \theta \in [-\pi, \pi].
\end{equation}

\vspace{0.6cm}

\noindent Proof: Necessity: Since $g(s) + f(s) \in H$ and $|| \frac{ \frac{1}{\gamma} f(s)}{f(s) + g(s)} ||_{\infty} < 1$, by Rouche's Theorem, we know that

\begin{equation}
[g(s) +f(s)] + \frac{1}{\gamma} e^{j \theta} f(s) \in H, \,\,\,\,\, \forall \theta \in [-\pi, \pi]
\end{equation}

Sufficiency: Now suppose on the contrary that $||S||_{\infty} \geq \gamma$, namely, $|| \frac{ \frac{1}{\gamma} f(s)}{f(s) + g(s)} ||_{\infty} \geq 1$. Since
$| \frac{ \frac{1}{\gamma} f(s)}{f(s) + g(s)} |_{s=j \omega} |$ is a contiunous function of $\omega$, and since

\begin{equation}
\lim_{\omega \rightarrow \infty} | \frac{ \frac{1}{\gamma} f(s)}{f(s) + g(s)} |_{s=j \omega} | = \frac{1}{\gamma} < 1
\end{equation}

\noindent there must exist $\omega_0$ such that

\begin{equation}
| \frac{ \frac{1}{\gamma} f(s)}{f(s) + g(s)} |_{s=j \omega_0} | = 1
\end{equation}

\noindent Therefore, there exists $\theta_0 \in [-\pi, \pi]$ such that

\begin{equation}
\{ g(s) +f(s) + \frac{1}{\gamma} e^{j \theta_0} f(s) \} |_{s = j \omega_0} = 0
\end{equation}

\noindent which contradicts the original hypothesis. This completes the proof.

\vspace{0.4cm}

\noindent {\bf Lemma 9}

For any $\delta \in (0, 1), \theta \in [-\pi, \pi]$, we have

\begin{equation}
W(s)=: \{g(s) + (1 + \delta e^{j \theta}) f(s) | g(s)\in K_g(s), f(s)\in K_f(s) \} \subset H \Longleftrightarrow
\end{equation}

\begin{equation}
J_{1111}, J_{1212}, J_{2222}, J_{2121}, J_{1112}, J_{1222}, J_{2221}, J_{2111}, J_{1211}, J_{2212}, J_{2122}, J_{1121} \in H
\end{equation}

\vspace{0.6cm}

\noindent Proof: Necessity is obvious. To prove sufficiency, note that $W(s)$ is a set of polynomials with complex coefficients, and with constant order $n$.
By Lemma 7, it suffices to show that

\begin{equation}
0 \not\in W(j \omega), \,\,\,\,\, \forall \omega \in R
\end{equation}

\noindent Since $0 \not\in W(j \omega_{\infty})$ for sufficiently large $\omega_{\infty}$, we only need to show that

\begin{equation}
0 \not\in \partial W(j \omega), \,\,\,\,\, \forall \omega \in R
\end{equation}

\noindent where $\partial W(j \omega)$ stands for the boundary of $W(j \omega)$ in the complex plane.

To construct $\partial W(j \omega)$, note that $\arg (1+ \delta e^{j \theta}) \in (- \frac{\pi}{2}, \frac{\pi}{2})$. Suppose now
$\omega \geq 0$ and $\arg (1+ \delta e^{j \theta}) \in [0, \frac{\pi}{2})$. Then by Lemma 6, we know that
$K_g (j \omega), K_f (j \omega)$ are rectangles with edges parallel to the coordinate axes. The four vertices of 
$K_g (j \omega)$ are $g_{11} (j \omega), g_{12} (j \omega), g_{21} (j \omega), g_{22} (j \omega)$, respectively;
and the four vertices of $K_f (j \omega)$ are $f_{11} (j \omega),$ $f_{12} (j \omega),$ $f_{21} (j \omega), f_{22} (j \omega)$, respectively.
$(1+ \delta e^{j \theta}) K_f (j \omega)$ is generated by rotating $K_f (j \omega)$ by 
$\arg (1+ \delta e^{j \theta})$ counterclockwisely, and then scaling by $|1+ \delta e^{j \theta}|$. Thus,
$W(j \omega) = K_g (j \omega) + (1+ \delta e^{j \theta}) K_f (j \omega)$ is a convex polygon with eight edges.
These edges are parallel to either the edges of $K_g (j \omega)$ or the edges of $(1+ \delta e^{j \theta}) K_f (j \omega)$.
Therefore, their orientations are fixed (independent of $\omega$). The eight vertices of $W(j \omega)$ are (clockwisely) $J_{1111} (j \omega)$, 
$J_{1112} (j \omega)$, $J_{1212} (j \omega)$, $J_{1222} (j \omega)$, $J_{2222} (j \omega)$, $J_{2221} (j \omega)$, $J_{2121} (j \omega)$, $J_{2111} (j \omega)$, respectively.

Now suppose on the contrary that there exists $\omega_0 \geq 0$ such that

\begin{equation}
0 \in \partial W(j \omega_0) 
\end{equation}

\noindent Without loss of generality, suppose

\begin{equation}
0 \in \{ \lambda J_{1111} (j \omega_0) + (1- \lambda) J_{1112} (j \omega_0) | \lambda \in [0, 1] \}
\end{equation}

\noindent Namely, there exists $\lambda_0 \in (0, 1)$ such that

\begin{equation}
\lambda_0 J_{1111} (j \omega_0) + (1- \lambda_0) J_{1112} (j \omega_0) = 0
\end{equation}

\noindent Since $J_{1111} (s), J_{1112} (s) \in H$, we have

\begin{equation}
\frac{d}{d \omega} \arg J_{1111} (j \omega) > 0, \,\,\,\,\,\,\,\, \frac{d}{d \omega} \arg J_{1112} (j \omega) > 0
\end{equation}

\noindent Thus\cite{Ran}

\begin{equation}
\hspace{-5cm} \frac{d}{d \omega} \arg [J_{1112} (j \omega) - J_{1111} (j \omega)] |_{\omega = \omega_0} =
\end{equation}

\begin{equation}
(1- \lambda_0) \frac{d}{d \omega} \arg J_{1111} (j \omega) |_{\omega = \omega_0} + \lambda_0 \frac{d}{d \omega} \arg J_{1112} (j \omega) |_{\omega = \omega_0} > 0
\end{equation}

\noindent This contradicts the fact that the edges of $W(j \omega)$ have fixed orientations. Thus

\begin{equation}
0 \not\in \partial W(j \omega) 
\end{equation}

\noindent Suppose now
$\omega \leq 0$ and $\arg (1+ \delta e^{j \theta}) \in (- \frac{\pi}{2}, 0]$. Then 
$K_g (j \omega), (1+ \delta e^{j \theta}) K_f (j \omega)$ are the mirror images (with respect to the real axis) of the corresponding sets
in the case of $\omega \geq 0$ and $\arg (1+ \delta e^{j \theta}) \in [0, \frac{\pi}{2})$. Therefore, following an identical line of arguments, we have

\begin{equation}
0 \not\in \partial W(j \omega) 
\end{equation}

\noindent The cases when $\omega \geq 0$ and $\arg (1+ \delta e^{j \theta}) \in (- \frac{\pi}{2}, 0]$ and when
$\omega \leq 0$ and $\arg (1+ \delta e^{j \theta}) \in [0, \frac{\pi}{2})$ are also symmetric with respect to the real axis.
Hence, we only need to consider the former case. In this case, $K_g (j \omega), K_f (j \omega)$ are rectangles with edges parallel to the coordinate axes.  
$(1+ \delta e^{j \theta}) K_f (j \omega)$ is generated by rotating $K_f (j \omega)$ by 
$| \arg (1+ \delta e^{j \theta}) |$ clockwisely, and then scaling by $|1+ \delta e^{j \theta}|$. Thus,
$W(j \omega) = K_g (j \omega) + (1+ \delta e^{j \theta}) K_f (j \omega)$ is a convex polygon with eight edges.
These edges are parallel to either the edges of $K_g (j \omega)$ or the edges of $(1+ \delta e^{j \theta}) K_f (j \omega)$.
Therefore, their orientations are fixed (independent of $\omega$). The eight vertices of $W(j \omega)$ are (clockwisely) $J_{1111} (j \omega)$, 
$J_{1211} (j \omega)$, $J_{1212} (j \omega)$, $J_{2212} (j \omega)$, $J_{2222} (j \omega)$, $J_{2122} (j \omega)$, $J_{2121} (j \omega)$, $J_{1121} (j \omega)$, respectively.
Thus, following a similar argument, we have

\begin{equation}
0 \not\in \partial W(j \omega) 
\end{equation}

\noindent This completes the proof.

The following theorem shows that, for an interval system, the maximal $H_{\infty}$ norm of its 
sensitivity function is achieved at twelve (out of sixteen) Kharitonov vertices.

\vspace{0.3cm}

\noindent {\bf Theorem 5}

Suppose $g_{ij} (s) + f_{ij} (s) \in H, \,\, i,j=1,2$. Then

\begin{equation}
\hspace{-5cm} \max \{ || \frac{f(s)}{f(s) + g(s)} ||_{\infty} | g(s)\in K_g(s), f(s)\in K_f(s) \} =
\end{equation}

\begin{equation}
\hspace{-4.5cm} \max \{ || \frac{f_{i_2 j_2} (s)}{f_{i_2 j_2} (s) + g_{i_1 j_1} (s)} ||_{\infty} | (i_1 j_1 i_2 j_2) = (1111), (1212),
\end{equation}

\begin{equation}
(2222), (2121), (1112), (1222), (2221), (2111), (1211), (2212), (2122), (1121) \}
\end{equation}

\vspace{0.6cm}

\noindent Proof: Since $g_{ij} (s) + f_{ij} (s) \in H, \,\, i,j=1,2$, by Kharitonov's Theorem\cite{Khar}, we know that 
$K_g(s) + K_f(s) \subset H$. Let

\begin{equation}
\hspace{-5cm} \gamma_1 = \max \{ || \frac{f(s)}{f(s) + g(s)} ||_{\infty} | g(s)\in K_g(s), f(s)\in K_f(s) \} 
\end{equation}

\begin{equation}
\hspace{-4cm} \gamma_2 = \max \{ || \frac{f_{i_2 j_2} (s)}{f_{i_2 j_2} (s) + g_{i_1 j_1} (s)} ||_{\infty} | (i_1 j_1 i_2 j_2) = (1111), (1212),
\end{equation}

\begin{equation}
(2222), (2121), (1112), (1222), (2221), (2111), (1211), (2212), (2122), (1121) \}
\end{equation}

\noindent Then apparently

\begin{equation}
\gamma_1 \geq \gamma_2 \geq 1
\end{equation}

Now suppose $\gamma_1 \neq \gamma_2$, namely, $\gamma_1 > \gamma_2$. Then there exists
$\gamma_0$ such that $\gamma_1 > \gamma_0 > \gamma_2$. Thus, for any 
$(i_1 j_1 i_2 j_2) \in \{(1111), (1212), (2222), (2121), (1112), (1222), (2221),$ $(2111), (1211), (2212), (2122), (1121) \}$, we have

\begin{equation}
|| \frac{f_{i_2 j_2} (s)}{f_{i_2 j_2} (s) + g_{i_1 j_1} (s)} ||_{\infty} < \gamma_0
\end{equation}

\noindent Hence, by Lemma 8, we have

\begin{equation}
g_{i_1 j_1} (s) + (1 + \frac{1}{\gamma_0} e^{j \theta}) f_{i_2 j_2} (s) \in H, \,\,\,\,\, \forall \theta \in [-\pi, \pi]
\end{equation}

\noindent By Lemma 9, we know that

\begin{equation}
\{ g(s) + (1 + \frac{1}{\gamma_0} e^{j \theta}) f(s) | g(s)\in K_g(s), f(s)\in K_f(s) \} \subset H, \,\,\,\,\, \forall \theta \in [-\pi, \pi]
\end{equation}

\noindent Therefore, by Lemma 8, for any $g(s)\in K_g(s), f(s)\in K_f(s)$, we have

\begin{equation}
|| \frac{f(s)}{f(s) + g(s)} ||_{\infty} < \gamma_0
\end{equation}

\noindent Namely

\begin{equation}
\max \{ || \frac{f(s)}{f(s) + g(s)} ||_{\infty} | g(s)\in K_g(s), f(s)\in K_f(s) \} < \gamma_0
\end{equation}

\noindent That is, $\gamma_1 < \gamma_0$, which contradicts $\gamma_1 > \gamma_0 > \gamma_2$. This completes the proof.

\vspace{6mm}
\section{Conclusions}

Some Kharitonov-like robust Hurwitz stability criteria have been established for a class of complex polynomial families with nonlinearly correlated perturbations. These results have been extended to the polynomial matrix case and non-interval D-stability case. Applications of these results in testing of robust strict positive realness of real and complex interval transfer function families have also been presented.


\begin{thebibliography}{9}


\bibitem{Khar}V.L.Kharitonov. Asymptotic stability of an equilibrium
  position of a family of systems of linear differential
  equations, Differential'nye Uravneniya, vol.14, 2086-2088, 1978.

\bibitem{Khar1}V.L.Kharitonov. The Routh-Hurwitz problem for families of polynomials and quasipolynomials, Izvetiy Akademii Nauk Kazakhskoi SSR, Seria fizikomatematicheskaia, vol.26, 69-79, 1979.

\bibitem{Holl}C.V.Hollot and R.Tempo. On the Nyquist envelope of an interval plant family, IEEE Trans. on Automatic Control, vol.39, 391-396, 1994.
   
\bibitem{Bart}A.C.Bartlett, C.V.Hollot and L.Huang. Root locations of an entire polytope of polynomials: It suffices to check the edges, Mathematics of Control, Signals, and Systems, vol.1, 61-71, 1988.
  
\bibitem{Fu}M.Fu and B.R.Barmish. Polytope of polynomials with zeros in a prescribed set, IEEE Trans. on Automatic Control, vol.34, 544-546, 1989.

\bibitem{wang1}L. Wang and L. Huang. Vertex results for uncertain systems, Int. J. Systems Science, vol.25, 541-549, 1994.

\bibitem{wang2}L. Wang and L. Huang. Extreme point results for strict positive realness of transfer function families, Systems Science and Mathematical Sciences, vol.7, 371-378, 1994.
 
\bibitem{Barmish}B. R. Barmish, C. V. Hollot, F. J. Kraus and R. Tempo. Extreme point results for robust stabilization of interval plants with first order compensators, IEEE Trans. on Automatic Control, vol.37, 707-714, 1992. 

\bibitem{Chap}H. Chapellat, M. Dahleh and S. P. Bhattacharyya. On robust nonlinear stability of interval control systems, IEEE Trans. on Automatic Control, vol.36, 59-67, 1991.

\bibitem{wang3}L. Wang and L. Huang. Finite verification of strict positive realness of interval rational functions, Chinese Science Bulletin, vol.36, 262-264, 1991.
 
\bibitem{Ack1}J.Ackermann. Uncertainty structures and robust stability analysis, Proc. of European Control Conference, 2318-2327,1991.

\bibitem{Ack2}J.Ackermann. Does it suffice to check a subset of multilinear parameters in robustness analysis? IEEE Trans. on Automatic Control, vol.37, 487-488, 1992.

\bibitem{das}S. Dasgupta. Kharitonov's theorem revisited, Systems and Control Letters, vol.11, No.4, 381-384, 1988.

\bibitem{Ran}A.Rantzer. Stability conditions for polytopes of polynomials, IEEE Trans. on Automatic Control, vol.37, 79-89, 1992.








\end{thebibliography}
\end{document}